\documentclass[11pt]{amsart}

\usepackage{RWWS}

\begin{document}
\title{Rozansky-Witten theory}
\author{Justin Roberts}
\address{Department of Mathematics, UC San Diego, 9500 Gilman Drive,
La Jolla CA 92093}
\email{justin@math.ucsd.edu}
\date{\today}

\begin{abstract}
Rozansky and Witten proposed in 1996 a family of new three-dimensional
topological quantum field theories, indexed by compact (or
asymptotically flat) hyperk\"ahler manifolds. As a byproduct they
proved that hyperk\"ahler manifolds also give rise to Vassiliev weight
systems. These may be thought of as invariants of hyperk\"ahler
manifolds, so the theory is of interest to geometers as well as to
low-dimensional topologists. This paper surveys the geometrical
construction of the weight systems, how they may be integrated into
the framework of Lie algebra weight systems (joint work with Simon
Willerton), their applications, and an approach to a rigorous
construction of the TQFTs (joint work with Justin Sawon and Simon
Willerton).
\end{abstract}

\maketitle


\section{Introduction}

The 1996 paper by Rozansky and Witten \cite{RW} opened up a very
interesting new area of geometry and topology. They wrote down a
physical path integral, a sigma model involving integration
over the space of all maps from a closed oriented $3$-manifold $M^3$
into a hyperk\"ahler manifold $X^{4n}$, which should give rise to a
numerical invariant $Z_X(M) \in \C$, at least in the case when $X$ is
compact or asymptotically flat. We can view this as a construction of
a family of invariants $Z_X(-)$ of $3$-manifolds, or of a family of
invariants $Z_{(-)}(M)$ of hyperk\"ahler manifolds. Part of the appeal
of the theory is the marriage of such different worlds.

Subsequent work by Kontsevich, Kapranov, Hitchin and Sawon
\cite{KontsRW, Kapranov, S,HS} has developed the theory from the
geometer's point of view. In particular, Kapranov showed that the
hyperk\"ahler condition is unnecessarily strong, and that the theory
works for {\em holomorphic symplectic manifolds}. This is the point of
view we take here.

There are three natural contexts for Rozansky-Witten theory, each with
its own motivations:

{\em 1. Geometry.} Perhaps the most concrete application of
Rozansky-Witten theory is to study hyperk\"ahler manifolds. Compact
hyperk\"ahler manifolds seem quite rare, though as yet there is no
guess at a classification theorem. The Rozansky-Witten invariants
amount to various complicated computations involving the curvature
tensor of the manifold, and as such can be related to characteristic
numbers. By combining the Rozansky-Witten invariants with the wheeling
theorem of Bar-Natan, Le and Thurston \cite{DylanPhD}, Hitchin and
Sawon obtained a formula for the $L^2$-norm of the curvature in terms
of a topological invariant, which is described in section 9. One
might hope to extract (perhaps using the TQFT structure) further
identities of this type, and obtain more constraints on the topology
of compact hyperk\"ahler manifolds.

A more extensive discussion of the geometrical point of view and its
potential applications, together with a list of problems, appears in
\cite{RS}; the present paper will concentrate on the other two viewpoints.

{\em 2. Vassiliev theory.} A compact holomorphic symplectic manifold
$(X^{4n}, \omega)$, together with a holomorphic vector bundle $E \ra
X$, defines a {\em Vassiliev weight system} $w_{X,E}: \A \ra \C$ on
the usual algebra $\A$ of Jacobi diagrams with an external circle.
This should be compared with the more familiar construction: a metric
(for example, semisimple) Lie algebra $\lie{g}$, together with a
representation $V$, defines a weight system $w_{\lie{g}, V}: \A \ra
\C$.

In sections 2-4 I will describe the differential-geometric
construction of the weight systems, following Kapranov. In sections
5-8 I describe how Simon Willerton and I have integrated these two
rather different-looking constructions, using the derived category of
coherent sheaves on $X$.

Our original motivation in studying Rozansky-Witten theory was in fact
to try to understand better the algebra $\A$ and its relatives. All
our intuition about such diagram algebras seems to be derived from the
{\em invariant theory of Lie algebras}, but as Vogel \cite{Vo} has
shown, there is more to them than that. In Rozansky-Witten theory,
diagrams behave more like {\em cohomology classes}, and this seems a
valuable alternative point of view. Whether we will gain true insight
remains to be seen.

{\em 3. TQFT.} A compact holomorphic symplectic manifold $(X^{4n},
\omega)$ defines a {\em $3$-dimensional topological quantum field
theory} $Z_X$, a functor from the category of $3$-dimensional
cobordisms to the category of finite-dimensional complex vector
spaces. Again, there is a more familiar example: the {\em Chern-Simons
theory} of Witten, Reshetikhin and Turaev \cite{WittenQFTJP,
ReshTuraev}, in which a semisimple Lie algebra $\lie{g}$ and a root of
unity $q=e^{2 \pi i/r}$ define a TQFT $Z_{\lie{g}, q}$. The work of
Rozansky and Witten does not actually give a rigorous construction of
such a theory, but Justin Sawon, Simon Willerton and I are currently
working on providing one: this is described in sections 10 and 11.

This theory is a really non-trivial new kind of TQFT in three
dimensions, and exhibits important differences from Chern-Simons
theory, most notably that it is {\em not semi-simple}, and therefore
does not satisfy the usual kinds of gluing and splitting
axioms. Though many generalisations of the original axioms of TQFT
have been investigated since their introduction twelve years ago,
there were until now no compelling examples to make such study
worthwhile.

The TQFT also appears to be related to quantum cohomology, quantum
$K$-theory, and mirror symmetry. Perhaps the techniques of {\em
$3$-dimensional} TQFT will prove useful in these areas, and this
theory will finally provide a link between the mysterious world of
quantum and Vassiliev invariants and the more concrete geometric world
of Gromov-Witten-type invariants, which seems to include almost all
the other modern topological invariants.


\section{Weight systems}

The {\em Kontsevich integral} is an invariant of framed oriented knots
in $S^3$ which takes values in (the completion of) a certain rational,
graded, algebra $\A$. This algebra is defined as the rational span of
vertex-oriented trivalent {\em Jacobi diagrams} which have a preferred
oriented circle, modulo the vertex-antisymmetry and IHX/STU
relations. The relations are like the skein relations in knot theory:
they relate diagrams which differ only locally. The pictures below
shows an example Jacobi diagram and the three relations. In planar
pictures such as these, vertex orientations are taken to be
anticlockwise. We grade diagrams by their total number of vertices;
this is {\em twice} the conventional grading.

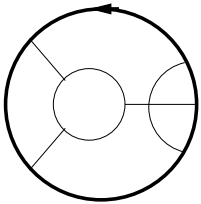
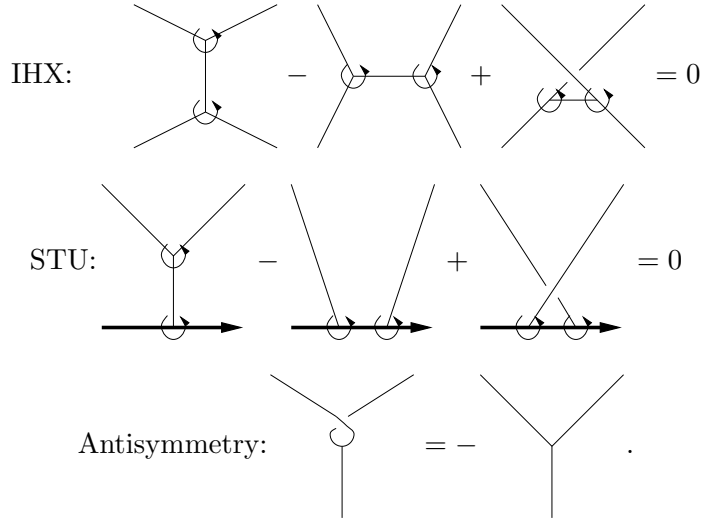
\[ \vpic{d} \]
\[ \text{IHX:}\qquad\vpic{I} - \vpic{H} + \vpic{X} =0\]
\[ \text{STU:}\vpic{S}-\vpic{T} +\vpic{U} =0\]
\[ \text{Antisymmetry:}\vpic{vxL} = - \vpic{vxR}.\]

A {\em weight system} is a linear map $\A \ra \Q$, or perhaps to some
other finite-dimensional vector space. Such maps enable us to study
the structure of $\A$, and by composing with the Kontsevich integral,
to construct more manageable scalar-valued knot invariants, the
Vassiliev invariants. This is described by Bar-Natan \cite{BN}.

The standard examples of weight systems, for a long time believed to
be essentially the only examples, are those coming from {\em metric
Lie algebras}. Let $\g$ be a finite-dimensional Lie algebra with a
non-degenerate invariant symmetric bilinear form $b$, and let $V$ be a
finite-dimensional representation of $\g$. These structures are
encoded by linear maps $\g \otimes \g \ra \g$, $\g \otimes \g \ra \Q$,
and $\g \otimes V \ra V$ respectively.  Using the metric to identify
$\g \cong \g^*$ enables us to transform these into equivalent but more
suitable structure tensors: a totally skew version $f \in
\Lambda^3\g^*$ of the bracket; the Casimir $c \in S^2\g$; and the
action $a_V \in V^* \otimes V \otimes \g^*$.

A Jacobi diagram of the kind spanning $\A$ defines a procedure for
contracting these tensors together and obtaining simply a rational
number. To do this, associate a copy of $f$ with each ``internal''
vertex, a copy of $a_V$ with each ``external'' vertex (those on the
preferred circle), and a copy of $c$ with each edge of the graph; now
contract the $\g$--$\g^*$ and $V$--$V^*$ pairs throughout. Note that
the orientation conventions correspond: we associate skew $3$-tensors
to skew vertices, and symmetric $2$-tensors to unoriented edges. The
Jacobi identity and the identity expressing that $V$ is a
representation of $\g$ correspond to the IHX and STU relations,
showing that this evaluation descends to the algebra $\A$ and gives a
weight system $w_{\g, V}: \A \ra \Q$.


\section{Chern-Weil theory on a complex manifold}

An attractive way to introduce Rozansky-Witten invariants is as a
generalisation of the usual Chern-Weil theory.  Suppose $E \ra X$ is a
smooth complex vector bundle on a smooth manifold $X$. Picking any
smooth connection on $E$ defines a curvature form $F \in \Omega^2(X;
\End(E))$, from which we seek to produce {\em topological} invariants
of $E$, quantities independent of the choice of connection. It is
well-known that the forms
\[\ch_d(E)= \frac{1}{d!}\tr \left\{\left(\frac{-F}{2 \pi i}\right)^d\right\}\] are closed, and that varying the connection
alters them by coboundaries, so that their de Rham cohomology classes
in $H^*(X;\C)$ are indeed invariants, namely the components of the
Chern character of $E$. The theory of characteristic classes shows
that these are the only functions of the curvature one need consider;
in any case, there are no other sensible ways to combine an
endomorphism-valued $2$-form with itself and end up with a
complex-valued differential form on $X$. A schematic picture of the
form representing $\ch_d(E)$ (in the case $d=8$)is shown below. Each
arrow represents a copy of $F$, the blob denoting the
``$2$-form-ness'' and the arrow shows that $F$, as a section of an
endomorphism bundle has one ``input'' and one ``output''. Eight
disjoint arrows would represent the wedge product $F^8 \in
\Omega^{16}(X; \End(E)^{\otimes 8})$; the circle illustrates the
contraction of consecutive inputs and outputs which computes the trace
of the product of the endomorphisms.
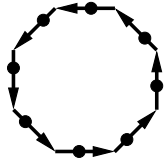
\[ \vpic{ch} \]

If we look at a bundle with a structure group other than $GL(n,\C)$
then there may be more subtle algebraic operations we can use; in
fact, the ring of invariant polynomials on the Lie algebra
parametrises the different (complex-valued) characteristic classes for
$G$-bundles. For example, in the case real oriented case with
structure group $SO(2n)$, the ring is generated by the traces of even
powers (corresponding to the Pontrjagin classes) together with the
Pfaffian, corresponding to the Euler class.

The basis of the Rozansky-Witten refinement (due, in this form, to
Kapranov) occurs when $E$ is a {\em holomorphic} bundle on a {\em
complex manifold} $X$. In this case there is a preferred class of
connections, namely those compatible with the holomorphic structure
and with some smooth hermitian metric on $E$. The curvature form $F_E$
of such a connection lies in $\Omega^{1,1}(\End(E))$, but can also be
thought of as an element $R_E \in \Omega^{0,1}(T^*\otimes\End(E))$,
where $T^*$ denotes the holomorphic cotangent bundle of $X$. Let $R_T$
denote a corresponding curvature form $R_T \in
\Omega^{0,1}(T^*\otimes\End(T))$ for the {\em holomorphic tangent bundle}
$T$. Each of these can be pictured as a {\em trivalent vertex} with
two input legs, one output leg, and a vertex carrying the
``$1$-form-ness''; larger trivalent graphs can then be used to index
the ways in which they can be combined tensorially, yielding a richer
range of possibilities than in the basic Chern-Weil case.
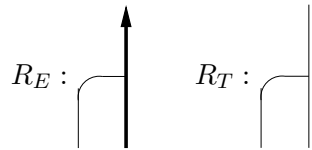
\[ R_E:\vpic{re}\qquad R_T:\vpic{rt} \]

In particular, the following three pictures denote elements of
$\Omega^{0,2}(T^*\otimes T^* \otimes \End(E))$ obtained by forming the
exterior products $R_T \wedge R_E, R_E \wedge R_E, R_E \wedge R_E$ and
then contracting tensorially according to the graph. The fundamental
lemma is that the sum of these three elements is a $\bar
\partial$-coboundary. Cohomologically, it will become the Jacobi or
IHX identity.
\[ \vpic{c1} + \vpic{c2} + \vpic{c3} = \bar\partial(\ldots).\]
A similar lemma shows that the element $R_T$ is symmetric in its two
inputs, up to a coboundary. Kapranov shows how this identity (which is
really just the Bianchi identity on $T^* \otimes E^* \otimes E$,
expanded using the Leibniz identity), defines the structure of an {\em
$L_\infty$-algebra} on the Dolbeault complex of forms with values in
$T$. We will not use this result here, preferring to work at the level
of cohomology.


\section{Rozansky-Witten weight systems}

Now let $(X^{4n}, \omega)$ be a {\em holomorphic symplectic manifold},
that is a complex manifold of {\em real} dimension $4n$, with $\omega
\in \Omega^0(\Lambda^2T^*)$ a holomorphic non-degenerate skew $2$-form on
$X$. Non-degeneracy of $\omega$ gives a holomorphic identification $T
\cong T^*$, and we can use this to convert the curvature form $R_T$
and symplectic form $\omega$ itself into alternative versions
\[ f_T \in \Omega^{0,1}(T^* \otimes T^* \otimes T^*),\qquad \omega^{-1} \in
\Omega^0(\Lambda^2T).\] The $1$-form $f_T$ turns out to be totally
symmetric in its three tensor indices.

Suppose we take a $v$-vertex trivalent graph with a preferred oriented
circle. By associating $f_T$ to its internal vertices, $\omega^{-1}$
to its edges, $R_E$ to its external vertices and wedging/contracting
accordingly, we obtain an element of $\Omega^{0, v}(X)$. This element
is $\bar \partial$-closed, and varying the connections on $E,T$ alters
it by a coboundary, giving a well-defined cohomology class. To avoid a
sign ambiguity here we need to orient the edges of the graph (because
$\omega^{-1}$ is skew) and order its vertices (because the vertices
carry $1$-forms, which anticommute). Remarkably, such an orientation
of a graph is naturally equivalent to a vertex-orientation of the
kind required in defining $\A$. The IHX and STU identities are
satisfied because of the Bianchi identity, and the result is that we
have a weight system
\[ w_{X,E} : \A^v \ra H_{\bar \partial}^{0,v}(X).\]

If $X^{4n}$ is {\em compact} then we can make a weight system $\A^{2n} \ra
\C$ by associating to a graph $\Gamma$ the integral
\[ b_\Gamma(X,E)=\int_X w_{X,E}(\Gamma) \wedge \omega^n.\]
These numbers are the original Rozansky-Witten invariants studied by
Sawon \cite{S}. Now compact holomorphic symplectic manifolds seem
quite rare: K\"ahler ones are hyperk\"ahler, by a theorem of Yau, and
the known examples comprise two infinite families (the Hilbert schemes
of points on the $4$-torus or $K3$ surface) and two exceptional
examples due to O'Grady; even non-K\"ahler ones are relatively few.

However, the construction of the cohomology-valued $w_{X,E}$ described
above does {\em not} require $X$ to be compact. Although in general we
may know little about the target Dolbeault cohomology groups, these
cases may turn out to be the most interesting for geometers, simply
because there are so many examples of (non-compact) moduli spaces
arising in geometry and physics which are naturally hyperk\"ahler, and
therefore holomorphic symplectic.


\section{Weight systems revisited}


What I will now explain is how to bring together the two extremely
different constructions of weight systems we have seen. The basic
ingredient is a reformulation of the construction of section 2 in a
way that generalises to metric Lie algebras {\em in categories other
than the category of vector spaces}. We will see subsequently that a
holomorphic symplectic manifold gives rise to such a category and Lie
algebra.

Suppose $\scrC$ is a {\em symmetric $\C$-linear tensor category}. This
is a category whose morphism sets are complex vector spaces, and which
is equipped with an associative tensor operation $\otimes: \scrC
\times \scrC \ra \scrC$, commutative in the sense that there are
isomorphisms $\tau_{A,B}: A\otimes B\to B\otimes A$, for all pairs of
objects $A,B$, satisfying $\tau_{A,B}\tau_{B,A} = \id$. The most
obvious example other than the usual category of complex vector spaces
is the category of {\em super}- (or {\em $\Z_2$-graded}) vector
spaces, in which the commutativity isomorphism incorporates a sign to
make odd elements {\em anticommute} with one another. It is the fact
that this commutativity isomorphism cannot be ignored that led Vogel
and Vaintrob \cite{Vo, Va} to the picture I am about to describe: for
Lie superalgebras, the standard construction explained earlier just
doesn't work.

An object $L$ in $\scrC$ is a {\em Lie algebra} if it is equipped with
a morphism $L \otimes L \ra L$ satisfying the Jacobi identity,
interpreted as a linear relation between three morphisms in
$\Mor(L^{\otimes 3}, L)$. Pictorially this is described by 
\[ \vpic{j1}+ \vpic{j2}+\vpic{j3}=0.\]
where the pictures are read from bottom to top, the trivalent vertex
represents the bracket, and the quadrivalent crossing represents the
commutativity morphism $\tau_{L,L}$. Such planar pictures are very
useful for describing compositions in a tensor category, and are very
common in the TQFT literature. See Bakalov and Kirillov \cite{BK}, for
example. 

Let $1$ be the unit object for the tensor product in $\scrC$. A {\em
non-degenerate metric} on a Lie algebra is a pair of morphisms $L
\otimes L \ra 1, 1 \ra L \otimes L$,  pictured using a ``cap'' and
``cup'', which satisfy a certain ``$S$-bend'' relation. A {\em left
module} $M$ over $L$ is an object with a morphism $L \otimes M \ra M$,
pictured as a trivalent vertex with a special oriented line marking
the legs corresponding to $M$.

A pair $(L,M)$ determines a weight system 
\[ w_{L,M}: \A \ra \Mor_\scrC(1,1). \]
To compute it, a diagram in $\A$ must first be drawn in the plane in a
Morse position with respect to the vertical coordinate. Slicing it
into horizontal bands, each containing a single vertex, cup or cap,
defines a way to compose the structural morphisms above, starting and
ending at the trivial object $1$. (An attempt at Morsifying the
example diagram from section 1 is shown below.) The result turns out
to be independent of the planar picture used (remember that the
diagrams in $\A$ are simply abstract trivalent graphs) and satisfies
the correct orientation and IHX relations.

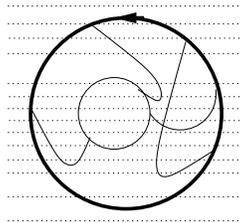
\[ \vpic{md}\]


\section{Derived categories}


The use of derived categories is necessary for the algebraic
unification we are proposing, and is even more important for the
construction of the TQFT. This section contains a brief description of
their structure, concentrating on the fact that {\em morphism sets in
a derived category are cohomology groups}. Because of this, the
derived world is the natural place to work if one wants to {\em
compose} lots of cohomology classes. For details, see Thomas
\cite{Th} and Gelfand and Manin \cite{GM}. 

Let $\scrC$ be an abelian category, for example the category of
modules over some ring $R$. Thus, $\scrC$ has direct sums, kernels and
cokernels, and it makes sense to consider the category of (bounded)
chain complexes $\Ch(\scrC)$. The most straightforward notion of
equivalence for chain complexes is {\em chain homotopy equivalence},
and it is straightforward to pass to the corresponding homotopy
category $K(\scrC)$. But this is not really the correct notion: better
is to use the {\em quasi-isomorphisms}, the chain maps which induce
isomorphisms on homology, to generate the equivalence. (Whitehead's
theorem in topology, that a map between simply-connected CW-complexes
inducing isomorphisms on homology is a homotopy-equivalence, is a
helpful justification here.) Now homotopy-equivalences are always
quasi-isomorphisms, but the converse is not true. For example, a
quasi-isomorphism from a complex of free modules to a complex of
torsion modules has no inverse, and this also shows that
quasi-isomorphism isn't an equivalence relation. Thus, by {\em
forcibly} symmetrising the quasi-isomorphisms in $K(\scrC)$, we obtain
a quotient category $D(\scrC)$ --- the derived category of $\scrC$.

It is characterised by the universal property that any functor defined
on $\Ch(\scrC)$ which takes quasi-isomorphisms to isomorphisms (the
most obvious example being the homology-group functors $h^i:
\Ch(\scrC) \ra \scrC$) factors through $D(\scrC)$. The {\em objects}
of the derived category are the same as those of $\Ch(\scrC)$. In
particular, objects of the original category $\scrC$ may be identified
with chain complexes whose only non-zero term lies in degree zero.

Under suitable conditions, a {\em functor} $F: \A \ra \B$ between
abelian categories may also be derived to a functor $D(F): D(\A) \ra
D(\B)$. The {\em classical derived functors} associated to $F$ are
just the composites of the homology functors $h^i$ with $D(F)$. For
example, if $\scrC$ is the category of $R$-modules and $\A$ the
category of abelian groups, then the ``$R$-invariant part'' functor
$\Gamma=\Hom_R(1,-): \scrC \ra \A$ has a derived functor denoted ${\bf
R}\Gamma$, whose classical derived functors $H^i=h^i{\bf R}\Gamma$ are
the standard cohomology functors for $R$-modules, and are computable
using resolutions. (In classical homological algebra there is an
inessential distinction between applying such a functor to an object
of $\A$, computing e.g. the cohomology of a module, and to an object
of $D(\A)$, computing e.g. the hypercohomology of a complex of
modules.)  A second example is that the bifunctor $\Hom_R(-,-):
\scrC^{\text{op}} \times \scrC \ra \A$ has a derived functor
$D(\scrC)^{\text{op}} \times D(\scrC) \ra D(\A)$ denoted ${\bf
R}\text{Hom}$ or $\Ext$, whose classical derived functors are the
bifunctors $\Ext_R^i(-,-)$.

The {\em morphisms} in $D(\scrC)$ are cohomology groups, in fact
\[ \Mor_{D(\scrC)} (A, B) =\Ext_{\scrC}^0(A,B). \]
Observing that for objects of $\scrC$ (rather than complexes)
$\Ext^0_{\scrC}(A,B) =\Hom_{\scrC}(A,B)$ shows that $\scrC$ is
embedded in $D(\scrC)$ as a full subcategory. More generally, we have
\[ \Mor_{D(\scrC)} (A, B[i]) =\Ext_{\scrC}^i(A,B).\]
Under this identification, the Yoneda product of cohomology classes
\[ \Ext^i(A,B) \otimes \Ext^j(B,C) \ra \Ext^{i+j}(A,C). \] becomes
{\em composition} of morphisms, shifted so as to make sense. The
construction of weight systems is, in effect, a complicated Yoneda
product, which is far better expressed as a composition of morphisms.


\section{Rozansky-Witten weight systems}


Let $X$ be a holomorphic symplectic manifold. We know that we can
construct weight systems $w_{X,E}$ out of holomorphic vector bundles
$E$ on $X$. In order to fit these into the framework of section 5, we
want to view these bundles as {\em modules} over a Lie algebra in some
category.

Ignore the symplectic structure for now and think of $X$ simply as a
complex manifold. Let $\OO_X$ be its {\em structure sheaf}, the sheaf
of germs of holomorphic functions on $X$. It is a sheaf of {\em
rings}, and the sheaves of germs of sections of vector bundles are
examples of sheaves of {\em $\OO_X$-modules}, in fact {\em
locally-free} ones.  If we consider the more general {\em coherent
sheaves of $\OO_X$-modules}, those which are locally quotients of
finite-rank locally-free sheaves, then we obtain an {\em abelian
category}. It is then a natural step to pass to the derived category,
which we denote by $D(X)$.

The objects of $D(X)$, then, are bounded chain complexes of coherent
sheaves, though thinking of them as complexes of holomorphic vector
bundles won't do too much harm. The morphism sets are cohomology
groups: in particular, for a pair of sheaves $E,F$ (viewed as
complexes concentrated in degree zero, as described earlier) we have
\[ \Mor_{D(X)}(E,F[i])=\Ext^i(E,F), \]
a convenient description of the usual sheaf cohomology $\Ext$-groups.
This language subsumes the earlier differential-geometric language: in
fact, when $E$ is a holomorphic vector bundle, we can identify the
Dolbeault cohomology groups, via sheaf cohomology groups, as morphism
sets:
\[ H^{0,i}_{\bar \partial}(X;E) = H^i(E) = \Ext^i(\OO_X, E) =
\Mor_{D(X)}(\OO_X, E[i]).\]
Notice that $\OO_X$ is the {\em unit object} in the tensor category
$D(X)$, and that the last term in this sequence may be thought of as a
kind of ``space of invariants'' of the object $E[i]$, by analogy with
the situation in a module category.

In order to construct a Lie algebra in this category we need a
sheaf-theoretic version of the curvature form of a bundle.  Suppose
$E$ is a holomorphic vector bundle on $X$. Its {\em Atiyah class}
$\alpha_E \in \Ext^1(E, E \otimes T^*)$ is the extension class
(obstruction to splitting) arising from the exact sequence
\[ 0 \ra E \otimes T^*  \ra JE \ra E \ra 0,\]
where $JE$ denotes the bundle of $1$-jets of sections of $E$. (See
\cite{A} for the details.) It
corresponds to the cohomology classes of the curvature forms $R_E$ and
$F_E$ used in section 3 under the isomorphisms
\[\Ext^1(E, E \otimes T^*) \cong H^1(E^*  \otimes E \otimes T^*) \cong
H^{0,1}_{\bar\partial}(E^* \otimes E \otimes T^*)\cong
H^{1,1}_{\bar\partial}(\End(E)).\]

For our purposes it is better to make the identification 
\[ \Ext^1(E, E \otimes T^*) \cong \Ext^1(E \otimes T, E) \cong
\Mor_{D(X)}(E \otimes T, E[1]),\]
because this begins to look like the structural map of a right
module. Applying the shift functor $[-1]$ to each side produces a
morphism
\[ \alpha_E: E \otimes T[-1] \ra E.\]
Repeating this for the tangent bundle $T$, with a shift of $[-2]$ gives
\[ \alpha_T: T[-1] \otimes T[-1] \ra T[-1].\]
This is in fact a Lie bracket: it is {\em skew} rather than symmetric,
because of the parity shift, and the earlier Bianchi identity becomes
the Jacobi identity. It is possible to extend the definition of the
Atiyah class $\alpha_E$ to arbitary coherent sheaves and even
complexes of them, so that we have the following theorem.

\begin{thm}
If $X$ is a complex manifold then the object $T[-1]$ is a Lie algebra
in the category $D(X)$, and all other objects in $D(X)$ are modules
over it.
\end{thm}

Note that this is true for {\em any complex manifold}. So what about
the symplectic structure? By viewing $\omega$ first as as a morphism
$\OO_X \ra \Lambda^2T^*$ and then dualising and shifting judiciously,
we construct a pair of morphisms
\[ \omega: T[-1] \otimes T[-1] \ra \OO_X[-2], \qquad \omega^{-1}:
\OO_X \ra T[-1] \otimes T[-1] [2].\]
These are actually {\em symmetric} rather than skew, again as a result
of the interaction of the shift with the notion of commutativity.

Because these morphisms have non-trivial ``degrees'', namely the
outstanding shifts $[\pm 2]$, we have to alter the category $D(X)$ in
order to be able to interpret them correctly as defining a
non-degenerate metric. Define the {\em graded derived category}
$\tilde D(X)$ to have the same set of objects as $D(X)$, but with the
space of morphisms $A \ra B$ being the graded vector space
$\Ext^*(A,B)$, instead of just $\Ext^0(A,B)$. Composition of morphisms
in $\tilde D(X)$ is graded bilinear.) After this ``deformation'', we
can view $\omega$ as a genuine morphism $T[-1]
\otimes T[-1] \ra \OO_X$, and hence as a metric.

This deformation seems pleasing rather than puzzling if one compares
it with the deformation quantization of the category of
representations of a semisimple Lie algebra arising from the
Knizhnik-Zamolodchikov equation, which is explained further in section
10. In summary:

\begin{thm}
If $X$ is a holomorphic symplectic manifold then $T[-1]$ is a metric
Lie algebra in the graded derived category $\tilde D(X)$, and all
other objects are modules over it.
\end{thm}

From this theorem and the methods of section 5 we get weight systems 
$w_{X,E}: \A \ra H^*(\OO_X)$ for each object $E$ of $D(X)$.


\section{Wheels and wheeling}


Let $\g$ be a complex semisimple Lie algebra, $U(\g)$ its universal
enveloping algebra, and $S(\g)$ its symmetric algebra. The {\em
Poincar\'e-Birkhoff-Witt theorem} says that the map $PBW: S(\g) \ra
U(\g)$, given by including $S(\g)$ into the tensor algebra $T(\g)$ and
then projecting, is an isomorphism of (filtered) vector spaces. It is
equivariant, so gives an isomorphism $S(\g)^\g
\cong U(\g)^\g =Z(\g)$. 

Now $S(\g)$, $U(\g)$ are algebras, and though they aren't isomorphic
as algebras (one is commutative and one non-commutative!), their
invariant parts {\em are} isomorphic as algebras, though {\em not} by
the restriction of the PBW map. The isomorphism is called the {\em
Duflo isomorphism}, and is built from the function
\[ j^{\frac12}(x) = {\det}^{\frac12} \left\{\frac{\sinh(\ad x/2)}{(\ad
x/2)}\right\}\] on the Lie algebra $\g$. One can define a map $S(\g)
\ra S(\g)$ by viewing this function as lying in the completion of
$S\g^*$ and applying the contraction $S(\g)^* \otimes S(\g) \ra
S(\g)$, which is a kind of convolution or ``cap product''. The Duflo
theorem is that the composite
\[ S(\g)^\g \stackrel{j^{\frac12}\cap}{\longrightarrow} S(\g)^\g
\stackrel{PBW}{\longrightarrow} U(\g)^\g\] is an isomorphism.

Bar-Natan, Le and Thurston \cite{DylanPhD} discovered a purely
diagrammatic version of this statement. Let $\B$ be the space spanned
by Jacobi diagrams with both trivalent and unitrivalent vertices,
subject to the usual orientation convention and relations. It is a
bigraded rational commutative algebra under disjoint union: write
$\B^{v,l}$ for the part with $v$ internal trivalent vertices, and $l$
legs.

This algebra $\B$ plays the role of $S(\g)^\g$; the original algebra
$\A$, with connect-sum of diagrams as its product, plays the role of
$U(\g)^\g$. The analogue of the PBW map is $\chi\colon \B \ra
\A$, defined by  sending an $l$-legged diagram to the average of the
$l!$ diagrams obtained by attaching its legs in all possible orders to
an oriented circle. This is an isomorphism of rational graded vector
spaces \cite{BN}, but not of algebras.

The analogue of the function $j^{\frac12}$ is the special {\em wheels
element}
\[ \Omega =\exp \left\{\sum_{i=1}^\infty b_{2i}w_{2i}\right\},\]
living in the completion of $\B$. Here, $w_{2i}$ is a {\em wheel}
diagram (a circular hub with $2i$ legs) and the $b_{2i}$ are versions
of Bernoulli numbers such that replacing $w_{2i}$ by the function
$\tr(\ad x)^{2i}$ recovers $j^{\frac12}(x)$. {\em Contraction} (we
will again write a cap product) with a diagram $C \in \B$ is the
operation $C
\cap: \B \ra \B$ given by summing over all attachments of legs of $C$
to legs of the target diagram. This extends to an action of the
completion of $\B$. The {\em wheeling theorem} is then that
\[ \B \stackrel{\Omega\cap}{\longrightarrow} \B \stackrel{\chi}{\longrightarrow} \A \]
is an isomorphism of algebras. 

A {\em metric} Lie algebra $\g$ defines weight systems $\B \ra
S(\g)^\g$ and $\A
\ra U(\g)^\g$, setting up a commuting diagram which intertwines the
wheeling and Duflo isomorphisms (see \cite{BGRT} for the details). (It
seems peculiar here that the Duflo isomorphism itself is true for {\em
all} Lie algebras, not just the metric ones.)

This whole picture has an analogue for holomorphic symplectic
manifolds. We construct in \cite{RbW} a pair of objects $S, U$ in
$\tilde D(X)$ which are the symmetric and universal enveloping
algebras of the Lie algebra $T[-1]$, in the sense that they are
associative algebras with appropriate universal properties. The
description of $S$ is straightforward, because one can take tensor
powers and form symmetrisers inside $\tilde D(X)$. Constructing $U$ is
not so straightforward, because $\tilde D(X)$ is not abelian and the
usual quotient construction therefore doesn't make sense. As one might
expect there is also a PBW isomorphism $S \cong U$ which does not
respect their structure as associative algebras.

The ``invariant parts'' of these two algebras are the algebras
$\Mor(\OO_X, S)$ and $\Mor(\OO_X, U)$, which are just the total
cohomology groups of these objects and are finite-dimensional graded
complex algebras. They appear in Kontsevich's remarkable paper
\cite{KontsDefQuant1} under the names $\HT^*$ (cohomology of
polyvector fields) and $\HH^*$ (Hochschild cohomology). To be precise,
we have
\[ \Mor( \OO_X, S) =\bigoplus_n \HT^n(X) = \bigoplus_n
\bigoplus_{i+j=n} H^i(\Lambda^jT)\]
\[ \Mor( \OO_X, U) =\bigoplus_n \HH^n(X) = \bigoplus_n \Ext_{X \times
X}^n(\OO_\Delta, \OO_\Delta)\] where $\OO_\Delta$ is the structure
sheaf of the diagonal in $X \times X$. The algebra structures here are
the natural ones: wedge product and Yoneda product, respectively. The
usual Dolbeault cohomology $H^*(\Lambda^*T^*)$ acts by contraction on
$\HT^*$, and Kontsevich proved that cap product with the {\em
root-A-hat-class}, the characteristic class corresponding to the power
series
\[ \left\{\frac{\sinh(x/2)}{x/2}\right\}^{\frac12} \]
gives a Duflo-style isomorphism of algebras
\[ \HT^* \stackrel{\hat A^{\frac12}\cap}{\longrightarrow} \HT^*
\stackrel{PBW}{\longrightarrow} \HH^*.  \]
A {\em holomorphic symplectic} $X$ supplies weight system maps $\B \ra
\HT^*$ and $\A \ra \HH^*$ intertwining the wheeling and Kontsevich
isomorphisms. (As in the Duflo case though, Kontsevich's isomorphism is
true for {\em any} complex manifold, and not just for holomorphic
symplectic ones.)


\section{The Hitchin-Sawon theorem}


\begin{thm}
Let $X^{4n}$ be a compact hyperk\"ahler manifold. Then there is an
identity relating the root-A-hat-genus of $X$, its volume and the
$L^2$-norm $\norm R \norm$ of its Riemann curvature tensor:
\[ {\hat A}^{\frac12}[X] = \frac{1}{(192\pi^2 n)^n} \frac{\norm R
\norm^{2n}}{\vol(X)^{n-1}}.\]
\end{thm}
This theorem, proved in \cite{HS}, is a striking identity between a
topological and a geometric invariant of $X$. Its proof is a
consequence of the wheeling theorem, and shows the potential of
Rozansky-Witten theory in geometry.

The basic idea of the proof is worth explaining. We use the
Rozansky-Witten weight systems in the form $w_X: \B^{v,l} \ra
H^{l,v}_{\bar \partial}(X)$. There are two numerical evaluations one
can consider: if $\Gamma \in \B^{2n,2n}$ is a diagram with $2n$
internal vertices and $2n$ legs then we can define 
\[ c_X(\Gamma) = \int_X w_X(\Gamma), \]
whereas if $\Gamma \in \B^{2n,0}$ is a ``closed'' trivalent diagram
with $2n$ vertices then we can define 
\[ b_X(\Gamma) = \int_X w_X(\Gamma) \wedge \omega^n.\]
The key observation is that $c_X(\Gamma) = b_X(\text{cl}(\Gamma))$ (up
to normalisations and signs) when $\text{cl}(\Gamma)$ is the {\em
closure} of $\Gamma$, the element obtained by summing over all ways of
pairing up the legs of $\Gamma$. (Remember that we are always really
working with linear combinations of diagrams.) This is a
straightforward consequence of the shuffle product formula for the top
power $\omega^n$ of a symplectic form evaluated on a list of $2n$
vectors.

If $\Gamma$ is chosen to be a {\em polywheel} diagram, that is, a
disjoint union of wheels with $2n$ legs in total, then $w_X(\Gamma)$
turns out to be a product of Chern character terms with total degree
$2n$, and so $c_X(\Gamma)$ is a linear combination of Chern numbers.
In particular, for the special wheels element $\Omega$ we get the
root-A-hat polynomial $w_X(\Omega) = \hat A^{\frac12}(TX)$, and hence
$c_X(\Omega) = \hat A^{\frac12}[X]$.

On the other hand, the closure $\text{cl}(\Omega)$ may be computed
using Jacobi diagram techniques worked out by Bar-Natan, Le and
Thurston
\cite{DylanPhD}, and it turns out to be $\exp(\Theta/24)$, where
$\Theta$ is the two-vertex theta-graph. A standard result in
hyperk\"ahler geometry expresses this graph in terms of the norm of
the curvature and gives the result. 

Huybrechts used this identity to prove a finiteness theorem for
hyperk\"ahler manifolds \cite{Huy}. It is to be hoped that further
identities of a similar nature will emerge from the TQFT structure of
RW theory, and help in the programme of classification of
hyperk\"ahler manifolds.


\section{Link invariants}


Let us turn now to the {\em topological} invariants arising from
Rozansky-Witten theory. Recall that a semisimple Lie algebra $\g$,
together with a representation $V$, defines a weight system $w_{\g,
V}: \A \ra \C$, and that composing this with the Kontsevich integral
should give us an invariant of framed knots in $S^3$. Because the
Kontsevich integral actually lies in the {\em completion} of $\A$, we
must first introduce a variable $\hbar$ of degree two and {\em
complete} the weight system into a graded map
\[ w_{\g, V}: \hat\A \ra \C[[\hbar]], \]
from which we obtain a $\C[[\hbar]]$-valued invariant $Z_{\g, V}$ of
framed knots in $S^3$. 

There is an alternative formulation of this invariant using the
framework of ribbon categories, invented by Turaev \cite{Tu} to handle
the invariants such as the Jones polynomial which arise from quantum
groups. A {\em ribbon category} is a braided tensor category with a
compatible notion of duality. It gives rise to representations of the
category of framed ``coloured'' tangles: one assigns an object of the
category (a ``colour'') to each string of a tangle, and then composes
the elementary structural morphisms in the category (braidings and
dualities) according to a Morse-theoretic slicing of the tangle into
crossings, cups and caps. In particular, one gets invariants of framed
coloured links, with values in the endomorphisms of the unit object of
the category.

The {\em Knizhnik-Zamolodchikov equation} builds an interesting ribbon
category from a semisimple Lie algebra. We start with the usual
symmetric tensor category of modules over $\g$, tensor it with
$\C[[\hbar]]$, and then use the monodromy of the KZ equation to
introduce a new, braided tensor (and in fact ribbon) structure. The
resulting category has non-trivial associator morphisms (it is
sometimes called a {\em quasitensor category}), but turns out to be
equivalent to the strictly associative category of representations of
the quantum group $U_\hbar (\g)$; see Bakalov and Kirillov \cite{BK},
for example. Applying Turaev's machinery to this category gives back
the $\C[[\hbar]]$-valued invariant $Z_{\g, V}$ of framed knots.

It turns out that there is an analogous picture for a holomorphic
symplectic manifold $X$. As we have seen, a vector bundle or sheaf $E$
on $X$ gives a weight system $w_{X,E}: \A \ra H^*(\OO_X)$, and
therefore a $H^*(\OO_X)$-valued Vassiliev invariant $Z_{X,E}$ of
framed knots.

\begin{thm} 
The graded derived category $\tilde D(X)$ of a holomorphic symplectic
manifold $X$ may be given a ribbon structure so that the
associated invariant of a knot coloured by $E \in D(X)$ is $Z_{X,E}$.
\end{thm}

The starting point of this construction is the usual derived category
$D(X)$, which is a symmetric tensor category. The analogue of
tensoring with $\C[[\hbar]]$ is its replacement by the graded version
$\tilde D(X)$, in which the shift $[2]$ plays the role of $\hbar$. The
$[2]$ is attached to the symplectic form $\omega$, just as in the Lie
algebra case, $\hbar$ may be thought of as attached to the metric. We
then use the Kontsevich integral, which underlies the KZ equation, to
define the ribbon structure $\tilde D(X)$.

By construction, the tensor structure on $\tilde D(X)$ is not strictly
associative. We don't know whether it is possible to make a gauge
transformation, as Drinfeld does in the case of quantum groups
\cite{DrinfeldQHA}, to a form which is strictly associative and has a
local, but more complicated, braiding. Drinfeld's transformation is
defined purely algebraically, and it is not clear how to derive it
from geometry in the way we would like. Equally remarkably, the
$\C[[\hbar]]$-valued invariants coming from Lie algebras turn out to
be the expansions, on setting $q=e^\hbar$, of Laurent polynomials in
$q$. Is there some analogous hidden structure to the Rozansky-Witten
link invariants?


\section{TQFT}


The basis of Rozansky and Witten's work is the path integral which
defines the partition function for a closed oriented $3$-manifold $M$:
\[ Z_X(M)= \int_{{\rm Map}(M,X)} e^{iS} d\phi.\]
The integral is over all smooth maps $\phi: M \ra X$, and the action
$S$ is an expression involving a Riemannian metric on $M$ and the
hyperk\"ahler metric on $X$. It turns out to be independent of the
metric on $M$, meaning that the associated quantum field theory, too,
is {\em topological}. The link invariants described above may also be
formulated using the path integral, by insertion of suitable ``Wilson
loop'' observables. 

We would like to have a rigorous construction of this TQFT, with which
we can calculate and explore. That is, we would like to build a TQFT
which shares the predicted properties of the genuine physical theory;
a ``Reshetikhin-Turaev invariant'' for our ``Witten invariant''. This
is underway in joint work with Justin Sawon and Simon Willerton.

To do it we have to work combinatorially. We will build the TQFT by
specifying its values on elementary pieces of surfaces and
$3$-manifolds, calculate for larger manifolds by means of
decompositions into such pieces, and prove that our specified
elementary data satisfies the coherence relations which ensure that
different decompositions of a $3$-manifold compute the same
invariant. This is not, in all honesty, a very satisfactory method,
but there's little reasonable alternative.

A standard TQFT, such as Witten's Chern-Simons theory
\cite{WittenQFTJP}, is a functor from the category of
$(2+1)$-dimensional cobordisms to the category of vector spaces. It is
typically a {\em tensor functor}, taking disjoint unions to tensor
product; and a {\em unitary} one, taking surfaces to hermitian vector
spaces. There is a way to enhance such a TQFT into an {\em extended}
or {\em $(1+1+1)$-dimensional} theory, in which we also assign a {\em
category} to each closed $1$-manifold. In the usual theories this
category is {\em semisimple}, with finitely many simple objects, and
is promptly replaced by a set of {\em colours} corresponding to these
objects. In Chern-Simons theory, the category for a single circle is
the (truncated) representation category of a quantum group at a root
of unity. Its simple objects correspond to the irreps of the classical
group which lie inside a certain compact Weyl alcove.

The Rozansky-Witten TQFT is somewhat different, and requires a
different formalism, which has been studied by Freed, Segal, Tillmann
and Khovanov \cite{Ti, Kh}, amongst others.  In this formalism we view
a $(1+1+1)$-dimensional TQFT as a functor from the $2$-category of
$3$-cobordisms with corners, to the $2$-category of linear
categories. Thus, a $1$-manifold is sent to a category; a ``vertical''
$2$-dimensional cobordism between $1$-manifolds defines a functor
between these categories; and a ``horizontal'' $3$-dimensional
cobordism between two surfaces which have common lower and upper
$1$-manifolds defines a natural transformation of functors. There's no
good reason here for imposing a tensor product axiom, unitarity, or
semisimplicity of the categories associated to $1$-manifolds.

Indeed, our theory has none of these. Its basic feature is that a
$1$-manifold consisting of $k$ circles gets sent to the graded derived
category $\tilde D(X^k)$: clearly $\tilde D(X)$ is not semisimple and
it is {\em not} true that $\tilde D(X \times X) \cong \tilde D(X)
\times \tilde D(X)$.

To construct the functor associated to a surface, we break it using
Morse theory into elementary pieces, and write down the basic functors
associated to the different types of handles. For example, a
$2$-handle gives the pushforward functor $\tilde D(X) \ra \tilde
D(\pt)$; a $0$-handle, the functor $\tilde D(\pt) \ra \tilde D(X)$
which sends $\C$ to $\OO_X$; and an index $1$-handle which joins two
circles, the functor $\tilde D(X \times X) \ra \tilde D(X)$ obtained
by taking derived tensor product with the structure sheaf of the
diagonal, then pushing forward once.

To prove independence of the decomposition we have to use the ribbon
structure on $\tilde D(X)$, and some extra properties resembling
modularity which arise from diagrammatic identities and the Kontsevich
integral. These are enough to check that the Moore-Seiberg equations
are satisfied. To define the invariants for elementary $3$-manifolds
is more complicated, and uses Walker's framework
\cite{Wa}.

It is possible to recover a traditional $(2+1)$-dimensional TQFT
functor which sends surfaces to vector spaces instead of functors from
our axioms. The {\em empty} $1$-manifold is sent to $\tilde D(\pt)$,
which is essentially the category of graded vector spaces. A {\em
closed} surface defines a functor $\tilde D(\pt)
\ra \tilde D(\pt)$, which when applied to the generating object $\C$ in
$\tilde D(\pt)$ outputs the desired graded vector space.

Using the basic data specified above, we can compute that the graded
vector space associated to a closed surface $\Sigma_g$ is isomorphic
(non-canonically) to the cohomology $H^*((\Lambda^*T)^{\otimes g})$,
as postulated by Rozansky and Witten. An action of the mapping class
group on this space emerges from the Moore-Seiberg equations.  The
superdimensions of these spaces are given for the sphere and torus by
the Todd genus and Euler characteristic of $X$, respectively, and they
vanish for $g \geq 2$. The space associated to $S^2$, which is
naturally an associative algebra in any TQFT, is in fact the ring
$H^*(\OO_X)$, which for an irreducible hyperk\"ahler $X$ is the
truncated polynomial ring $\C[\bar \omega]/(\bar \omega^{n+1})$. It is
not semisimple, which shows that the TQFT cannot be unitary.

This TQFT has a close relationship with the Le-Murakami-Ohtsuki
invariant \cite{LMO}, which is a kind of extension of the Kontsevich
integral to $3$-manifolds, with values in the algebra of trivalent
Jacobi diagrams. In principle, the TQFT invariant of a closed
$3$-manifold should equal the LMO invariant, evaluated using a weight
system coming from $X$; the TQFT itself should be the evaluation of
the universal TQFT constructed by Murakami and Ohtsuki \cite{MO}, which
underlies the LMO invariant. However, their TQFT is somewhat badly
behaved axiomatically, and there are subtle differences in
normalisation which make it unsatisfactory to try to use this approach
as a definition. The $(1+1+1)$-dimensional framework seems to work so
well when we are building the specific theory associated to $X$ that
it is better to use it as a foundation, and to make connections with
the universal theory in retrospect.


\section{Future directions}


In \cite{RS} we give an extensive discussion of the potential
applications of RW invariants to the geometry of hyperk\"ahler (and
certain other types of) manifolds, and conclude with a problem
list. Further discussion of the interaction with Vassiliev theory, the
theory of quantum invariants and TQFT will appear in \cite{RSW}.  But
it is worth mentioning a few problems here, just to give the flavour
of possible future research.

{\em What is the meaning of the $\hat A^{\frac12}$ genus?}

This genus appears in Kontsevich's theorem and the Hitchin-Sawon
theorem, as we have seen. But it also appears in physics, in the
context of cohomological $D$-brane charge. It seems natural that these
two occurrences are related, but how and why? (With impeccable timing,
physicists are also now adopting derived categories of coherent
sheaves as basic structures in string theory!) A straightforward
question is whether the root-A-hat genus, like the honest A-hat genus,
has interesting integrality properties, and can be interpreted for
manifolds with some appropriate geometric structure as the index of a
natural elliptic operator.

{\em Are the Rozansky-Witten weight systems new, and can their
associated Vassiliev invariants detect orientation of knots?}

Vogel showed that the primary examples of Vassiliev weight systems,
those coming from complex semisimple Lie algebras and superalgebras,
do not span the whole space of weight systems.  But we do not yet know
whether the Rozansky-Witten weight systems lie outside their span or
not.

It's easy to prove that the Vassiliev invariants coming from Lie
algebra weight systems are unable to distinguish knots from their
reverses. It is in fact thought likely that {\em no} Vassiliev
invariants can separate knots from their reverses, and there is an
alternative purely diagrammatic statement of this conjecture. But the
proof that Lie algebra weight systems fail doesn't work for
Rozansky-Witten invariants, so there is still potential here.

{\em Is there a geometric quantization approach to the vector spaces
$Z(\Sigma_g)$?}

These vector spaces, as we have seen, can be constructed from a
combinatorial approach to the TQFT. This is a reasonable but crude
approach: much better would be to give a direct geometric (and
completely functorial) construction. In Witten's Chern-Simons theory
with gauge group $G$, there is an approach via geometric quantization
of the moduli space of flat $G$-connections on a surface. Is there an
analogue in Rozansky-Witten theory? A straightforward guess is that it
might be possible to define a virtual structure sheaf of the moduli
space of holomorphic maps from a closed Riemann surface $\Sigma_g$ to
$X$, and take its cohomology as the graded vector space. We would then
require a ``projectively flat connection'', a system of coherent
isomorphisms between these sheaves, as the complex structure on
$\Sigma_g$ varies. Going a level deeper, how does one construct the
category $\tilde D(X)$ functorially from a given circle and manifold
$X$? In the Chern-Simons case, one uses the representation category of
the group of loops from the circle to $G$. What is the analogue? We
would require surfaces with boundary to generate functors between such
categories by means of a Fourier-Mukai transform operation.


\noindent{\bf Acknowledgements} The research described here was
partially supported by an EPSRC Advanced Fellowship and NSF grant
DMS-0103922. I would like to thank Gordana Matic and the UGA
topologists for hosting such an excellent conference.




\begin{thebibliography}{BGRT}

\bibitem[A]{A} M.~F.~Atiyah. {\em Complex analytic connections in
        fibre bundles}, Trans.\ Amer.\ Math.\ Soc.\ 85 (1957),
        181-207.

\bibitem[BK]{BK} B.~Bakalov, A.~Kirillov. {\em Lectures on tensor
        categories and modular functors} University Lecture Series 21,
        AMS (2001).

\bibitem[BN]{BN} D.~Bar-Natan. {\em On the
        Vassiliev knot invariants}, Topology 34 (1995), 423-472.

\bibitem[BGRT]{BGRT} D.~Bar-Natan, S.~Garoufalidis, L.~Rozansky,
        D.~P.~Thurston. {\em Wheels, wheeling, and the Kontsevich
        integral of the unknot}, Israel J.\ Math.\ 119 (2000),
        217-237.

\bibitem[Dr]{DrinfeldQHA} V.~G.~Drinfel'd. {\em Quasi-Hopf
        algebras}, Leningrad J.\ Math.\ 1 (1990) 1419-1457.

\bibitem[GM]{GM} S.~I.~Gelfand, Yu.~I.~Manin. {\em Homological
        algebra}, Springer (1999).

\bibitem[HS]{HS} N.~Hitchin, J.~Sawon. {\em Curvature and
        characteristic numbers of hyper-K\"ahler manifolds}, Duke
        Math.\ J.\ 106 (2001), no. 3, 599-615.

\bibitem[H]{Huy} D.~Huybrechts. {\em Finiteness results for
	hyperk\"ahler manifolds}, arXiv preprint {\tt math.AG/0109024}

\bibitem[Ka]{Kapranov} M.~Kapranov. {\em Rozansky-Witten
        invariants via Atiyah classes}, Compositio Math.\ 115 (1999),
        71-113.

\bibitem[Kh]{Kh} M.~Khovanov. {\em A functor-valued invariant of
	tangles}, arXiv preprint {\tt math.QA/0103190}

\bibitem[Ko1]{KontsDefQuant1} M.~Kontsevich. {\em Deformation
        quantization of Poisson manifolds I}, arXiv preprint {\tt
        q-alg/9709040}.

\bibitem[Ko2]{KontsRW} M.~Kontsevich. {\em Rozansky-Witten invariants
        via formal geometry}, Compositio Math.\ 115 (1999), 115-127.

\bibitem[LMO]{LMO} T.~Q.~T.~Le, J.~Murakami, T.~Ohtsuki. {\em On a
        universal perturbative invariant of $3$-manifolds}, Topology
        37 (1998), 539-574.



\bibitem[MO]{MO} J.~Murakami, T.~Ohtsuki. {\em Topological quantum
        field theory for the universal quantum invariant}, Comm.\
        Math.\ Phys.\ 188 (1997), 501-520.

\bibitem[RT]{ReshTuraev} N.~Y.~Reshetikhin, V.~G.~Turaev. {\em
        Invariants of 3-manifolds via link polynomials and quantum
        groups}, Invent.\ Math.\ 103 (1991), 547-597.

\bibitem[RS]{RS} J.~D.~Roberts, J.~Sawon. {\em Generalisations of
	Rozansky-Witten invariants}, in preparation.

\bibitem[RSW]{RSW} J.~D.~Roberts, J.~Sawon, S.~Willerton. {\em On the
        Rozansky-Witten TQFT}, in preparation.

\bibitem[RW]{RbW} J.~D.~Roberts, S.~Willerton. {\em On the
	Rozansky-Witten weight systems}, in preparation.

\bibitem[RzW]{RW} L.~Rozansky, E.~Witten. {\em
        Hyper-K\"ahler geometry and invariants of $3$-manifolds},
        Selecta Math.\ {\bf 3} (1997) 401--458.

\bibitem[Sa]{S} J.~Sawon. {\em Rozansky-Witten invariants of
	hyperk\"ahler manifolds}, PhD Thesis, University of Cambridge,
	1999.

\bibitem[Tho]{Th} R.~P.~Thomas. {\em Derived categories for the working
        mathematician}, arXiv preprint {\tt math.AG/0001045}.

\bibitem[Thu]{DylanPhD} D.~P.~Thurston. {\em Wheeling: a diagrammatic
        analogue of the Duflo isomorphism}, UC Berkeley thesis (2000),
        arXiv preprint {\tt math.QA/0006083}.

\bibitem[Ti]{Ti} U.~Tillmann. {\em S-Structures for k-linear
	categories and the definition of a modular functor}, arXiv
	preprint {\tt math.GT/9802089}

\bibitem[Tu]{Tu} V.~G.~Turaev. {\em Quantum invariants of knots and
        3-manifolds}, de Gruyter Studies in Mathematics, 18 (1994).

\bibitem[Va]{Va} A.~Vaintrob. {\em Vassiliev knot invariants and
        Lie $S$-algebras.}, Math.\ Res.\ Lett.\ 1 (1994), 579-595.

\bibitem[Vo]{Vo} P.~Vogel. {\em Algebraic structures on modules of
        diagrams}, preprint (1995).

\bibitem[Wa]{Wa} K.~Walker. {\em On Witten's 3-manifold
        invariants}, preprint (1991).

\bibitem[Wi]{WittenQFTJP} E.~Witten. {\em Quantum field
        theory and the Jones polynomial}, Comm.\ Math.\ Phys.\ 121
        (1989), 351-399.

\end{thebibliography}
\end{document}